\documentclass[a4paper, 12pt]{article}
\usepackage{amssymb, amsmath, amsfonts}
\usepackage[english,russian]{babel}
\usepackage{amsthm,amscd,amsfonts,latexsym,color,epsfig}

\begin{document}

\numberwithin{equation}{section}

\begin{center}
RIEMANN--LIOUVILLE DERIVATIVE APPROACH OF THE FRACTIONAL X-RAY 
DIFFRACTION THEORY: THE CAUCHY NONLOCAL PROBLEM

\end{center}
\begin{center}
{\bf  Murat O. Mamchuev, Felix N. Chukhovskii}
\end{center}


\small

{\bf Abstract.} 
Towards the theory of the fractional X-ray diffraction crystal optics,  the Riemann--Liouville derivative approach  has been developed and analyzed for arbitrary fractional-order-parameter $\alpha,$ $\alpha\in (0, 1].$ Using the Green function 
technique,  fundamentals of the fractional matrix differential theory of the X-ray diffraction scattering by the imperfect crystals, have been formulated  in the terms of solving the Cauchy  nonlocal  problem.  
In a general case of an arbitrary displacement imperfect-crystal-lattice function $f(r)$, one has used the Liouville--Neumann-type series formalism to build up the resolvent-function-solution of the Cauchy problem rewritten in the form of the matrix integral Fredholm--Volterra equation of the second kind. 
A comparison was made between the theoretical results obtained using an approach based on the Riemann--Liouville derivative operator and the corresponding results obtained using fractional derivatives in the Gerasimov--Caputo sense.
In the case when the displacement field function $f(r)$ is the lateral coordinate $x$ is the $1D$ linear function,  $f(r)=ax,$  $a=const,$ one has obtained the rigorous Cauchy problem solutions expressed through the Mittag-Leffler-type and generalized trigonometric functions. 
One has numerically calculated, built up and analyzed in details,  the graph-patterns of the amplitudes of the transmitted in details d and diffracted waves along the crystal depth $t$ for different values of the fractional-order-parameter  $\alpha$ in the both cases of the theoretical  Riemann--Liouville derivative and Gerasimov--Caputo's one approaches.  Based on Mittag-Leffler-type asymptotic expressions (for large values of $t$), the decay behavior of the X-ray wave amplitude was determined for both the Riemann--Liouville and Gerasimov--Caputo derivative approaches.
 \medskip

{\bf MSC 2010\/}: Primary 35F35;
                  
                  Secondary 35F40, 35A08, 35C15, 35L40, 45F05, 35Q70, 35Q92;

 \smallskip

{\bf Key Words and Phrases}:  Diffraction optics system of fractional partial differential equations; the Riemann--Liouville fractional derivative; the  Cauchy nonlocal problem; the matrix resolvent-function-solution; the Fredholm--Volterra  integral equation; the generalized trigonometric functions; the generalized Euler's formula;

\normalsize

\section{Introduction}\label{sec1}

$\quad$ As is known, in recent years, to model and analyze the real-world phenomena in the terms of the mathematical physics, fractional calculus has become widespread of allowing to go beyond the integer differential equations solutions. In particular, the latter is important to taking into account and adequately considering such the effects as memory and nonlocality of the quantum radiation-medium interactions 
(see \cite{Nakhushev-2003-di}-\cite{Mamchuev-Chukhovskii-2023-fcaa} and huge cited references therein).

Recently \cite{Mamchuev-Chukhovskii-2023-fcaa}, 
\cite{Mamchuev-Chukhovskii-2025-jms}, the theory of the fractional X-ray diffraction crystal optics has been developed with using the  system of the fractional Takagi--Taupin (TT) differential equations in the terms of  the Gerasimov--Caputo (GC) derivatives of an arbitrary fractional order $\alpha,$ $\alpha\in (0, 1]$, along the crystal depth. 
In \cite{Mamchuev-Chukhovskii-2025-fcaa}, the integral formalism of the $3D$ fractional X-ray diffraction crystal optics has been done, and the matrix resolvent solution of the fractional GC-Cauchy problem has been constructed  as well.

Accordingly, in the case when the lattice-crystal displacement field function $f(r)$ vs the lateral coordinate x is the $1D$ linear function,  $f(r)=ax$,  $a=const$, the rigorous  analytical solution of the fractional in GC-Cauchy problem with an arbitrary fractional-order-parameter (FOP) $\alpha,$ $\alpha\in (0, 1]$  has been found out and analyzed.

As was pointed out in [7], various kinds of real problems-in-physics such as signal processing,  the radiation propagation,  atom diffusion in material media, are currently modeled by the aid of fractional derivative by highlighting in particular the Riemann--Liouville (RL) derivative and GC derivative approaches as the fundamental formulations in fractional calculus that generalizes the concept of derivatives to non-integer.

Based on the fractional TT differential equations, a goal of the present study is to push ahead the theory of the fractional X-ray diffraction crystal optics in the terms of the RL derivative concept beyond the frame of the GC derivative constructive scheme \cite{Mamchuev-Chukhovskii-2023-fcaa}, 
\cite{Mamchuev-Chukhovskii-2025-jms}.

In the paper, using the standard Green's function technique, a key of working is to develop the matrix Fredholm--Volterra integral equation formalism earlier proposed by authors 
\cite{Mamchuev-Chukhovskii-2023-fcaa}, 
\cite{Mamchuev-Chukhovskii-2025-jms}, and builds  up the matrix resolvent-function solution of the RL-Cauchy nonlocal problem. 

In the case of the lattice-crystal displacement field function $f (r)$ is the 1D linear function, namely: $f(r) = ax$, $a=const$, for an arbitrary  FOP $\alpha,$ $\alpha\in (0, 1]$ of the RL derivative, the exact analytical solution of the RL-Cauchy nonlocal problem has been found out and analyzed of comparing it with the GC concept-Cauchy problem solution. 

Thus, following up the theoretical formalism earlier developed in 
\cite{Mamchuev-Chukhovskii-2023-fcaa}, 
\cite{Mamchuev-Chukhovskii-2025-jms} the system of the original TT differential equations in the frame of the RL concept of the fractional X-ray diffraction crystal optics takes the form
$$\left (
\begin{array}{cc}
D_{0t}^{\alpha}-\frac{\partial}{\partial x}	& 0  	 \\
 0   	&      D_{0t}^{\alpha}+\frac{\partial}{\partial x}
\end{array}
\right )  
\left (
\begin{array}{c}
E_0(x,t)	 \\
E_h(x,t)
\end{array}
\right ) =$$
\begin{equation}\label{eq11}
=i
\left (
\begin{array}{cc}
\gamma						& \sigma \exp[if(x,t)]  	 \\
\sigma\exp[-if(x,t)]    	& \gamma	     
\end{array}
\right )  
\left (
\begin{array}{c}
E_0(x,t)	 \\
E_h(x,t)
\end{array}
\right ),
\end{equation}
with the Cauchy problem's condition
\begin{equation}\label{eq12}
\lim\limits_{t\to 0}D_{0t}^{\alpha-1}\left (
\begin{array}{c}
E_0(x,t)	 \\
E_h(x,t)
\end{array}
\right ) =
\left (
\begin{array}{c}
\varphi_0(x)	 \\
\varphi_h(x)
\end{array}
\right )\equiv \varphi(x), 
\quad -\infty<x<\infty,
\end{equation}
where $\varphi_0(x)$ and $\varphi_h(x)$ are the given real-valued functions.

The RL fractional integro-differential operator
of order $\alpha<0$ \cite{Nakhushev-2003-di} is equal to  
$$D_{at}^{\alpha}g(t)=\frac{\mathop{\rm sgn}(t-a)}{\Gamma(-\nu)}
\int\limits_a^t \frac{g(s)ds}{\lvert t-s\rvert^{\alpha+1}}, 
$$
and for $\alpha \in (0, 1]$ the operator $D_{at}^{\nu}$ is determined by 
following relation
$$D_{at}^{\alpha}g(t)=\mathop{\rm sgn}(t-a)\frac{d}{dt}D_{at}^{\alpha-1}g(t), 
$$
here $\Gamma(z)$ is the gamma-function.

\section{Integral fundamentals of the RL Cauchy non-local problem}\label{sec:3} 
\setcounter{equation}{0}

Let us convert the Cauchy non-local  problem in the 'differential form' (\ref{eq11})--(\ref{eq12}) to the matrix Fredholm--Volterra-type integral equation of the second kind. The latter takes a special sense to build up the resolvent-function solution of the Cauchy problem, which in turn is very important to computer-aided modelling
and recovering the crystal displacement field function $f({\bf R})$ from the microtomography data.

The system of differential TT-type equations (\ref{eq11}) may be rewritten into the form
\begin{equation}\label{eq31}
\left (
\begin{array}{cc}
O_{-}^{\alpha}-i\gamma		& 0  	 \\
 0   	&     O_{+}^{\alpha}-i\gamma	
\end{array}
\right ) {\bf E}=i\sigma {\bf KE},
\end{equation}
with initial condition
$$
\lim\limits_{t\to 0}D_{0t}^{\alpha-1}{\bf E}(x,t) =
\varphi(x), 
\quad -\infty<x<\infty,
$$
where
$$
O_{+}^{\alpha}=D_{0t}^{\alpha}+\frac{\partial}{\partial x},
\quad O_{-}^{\alpha}=D_{0t}^{\alpha}-\frac{\partial}{\partial x},
\quad
{\bf K}\equiv
{\bf K}(x,t)=
\left (
\begin{array}{cc}
0			  	&  e^{if(x,t)}  	 \\
e^{-if(x,t)}    & 0  
\end{array}
\right ),
$$
$$
{\bf E}\equiv
{\bf E}(x,t)=
\left (
\begin{array}{c}
E_0(x,t)  	 \\
E_h(x,t) 
\end{array}
\right ),
\quad
\varphi(x)=\left (
\begin{array}{c}
\varphi_0(x)	 \\
\varphi_h(x)
\end{array}
\right ). 
$$
Acting onto both sides of (\ref{eq31}) by the operator 
$diag(O_{+}^{\alpha}-i\gamma	, O_{-}^{\alpha}-i\gamma	),$ one obtains 
$$\left (
\begin{array}{cc}
O_{+}^{\alpha}-i\gamma		& 0  	 \\
 0   	&      O_{-}^{\alpha}-i\gamma	
\end{array}
\right )  
\left (
\begin{array}{cc}
O_{-}^{\alpha}-i\gamma		& 0  	 \\
 0   	&     O_{+}^{\alpha}-i\gamma	
\end{array}
\right )  
{\bf E}=$$
\begin{equation}\label{eq33}
=
i\sigma 
\left (
\begin{array}{cc}
O_{+}^{\alpha}-i\gamma		& 0  	 \\
 0   	&      O_{-}^{\alpha}-i\gamma	
\end{array}
\right )
{\bf KE}
\end{equation}
and after some straightforward routine calculations, one finds out
$$\left (\!\!
\begin{array}{cc}
O_{+}^{\alpha}-i\gamma		& 0  	 \\
 0   	&      O_{-}^{\alpha}-i\gamma	
\end{array}
\!\! \right )
{\bf KE}=
\left [D_{0t}^{\alpha}+
\Lambda
\frac{\partial}{\partial x}
-i\gamma \right ]({\bf KE})=
D_{0t}^{\alpha}({\bf KE})+$$
$$+
\Lambda
\left[
\left(\frac{\partial}{\partial x}{\bf K}\right){\bf E}+
{\bf K}\left(\frac{\partial}{\partial x}{\bf E}\right)
\right]
+{\bf K}\left(D_{0t}^{\alpha}{\bf E}\right)-
{\bf K}\left(D_{0t}^{\alpha}{\bf E}\right)
-i\gamma\,{\bf KE}=$$
$$=
{\bf K}
\left (\!\!
\begin{array}{cc}
O_{-}^{\alpha}-i\gamma	&     0 	 \\
0  	&     O_{+}^{\alpha}-i\gamma
\end{array}
\!\!\right )  
{\bf E}+D_{0t}^{\alpha}({\bf KE})+
\Lambda
\left(\frac{\partial}{\partial x}{\bf K}\right){\bf E}-
{\bf K}\left(D_{0t}^{\alpha}{\bf E}\right)=
$$
\begin{equation}\label{eq34}
={\bf K}
\left (
\begin{array}{cc}
O_{-}^{\alpha}-i\gamma	&     0 	 \\
0  	&     O_{+}^{\alpha}-i\gamma
\end{array}
\right ){\bf E}
+\left(D_{0t}^{\alpha}+if'_x\right)({\bf KE})
-{\bf K}\left(D_{0t}^{\alpha}{\bf E}\right),  
\end{equation}
where 
$\Lambda=\left (
\begin{array}{cc}
1	&    0   	 \\
0   &    -1  
\end{array}
\right ). $
Taking into account that the ${\bf E=E}(x,t)$ is the solution of Eq. (\ref{eq31}), the equality $\Lambda {\bf K}=-{\bf K}\Lambda$ and that  ${\bf K}^2$ is equal to unit matrix, from Eq. (\ref{eq34}) it directly follows
$$i\sigma\left (
\begin{array}{cc}
O_{+}^{\alpha}-i\gamma		& 0  	 \\
 0   	&      O_{-}^{\alpha}-i\gamma	
\end{array}
\right )
{\bf KE}=-\sigma^2 {\bf E}
+i\sigma\left(D_{0t}^{\alpha}+if'_x\right)({\bf KE})
-i\sigma{\bf K}\left(D_{0t}^{\alpha}{\bf E}\right).$$
The last allows us to rewrite down Eq. (\ref{eq33}) into the form
$$\left (
\begin{array}{cc}
(O_{+}^{\alpha}-i\gamma	)(O_{-}^{\alpha}-i\gamma	)+\sigma^2      & 0  	 \\
 0   	&      (O_{-}^{\alpha}-i\gamma	)(O_{+}^{\alpha}-i\gamma	)+\sigma^2
\end{array}
\right )  
{\bf E} =$$
\begin{equation}\label{eq35}
=i\sigma\left(D_{0t}^{\alpha}+if'_x\right)({\bf KE})
-i\sigma{\bf K}\left(D_{0t}^{\alpha}{\bf E}\right).  
\end{equation}

Further, we denote the Fourier transform of the function $f(x)$ by $(f(x))_k$,  the Laplace transform of the function $g(t)$ by $(g(t))_p$,  and  respectively, the double Fourier--Laplace transform of the function $h(x,t)$ by
$(h(x,t))_{k,p}$. 

By using the following formula
$$[D_{0t}^{\alpha}H(x,t)]_{p}=p^{\alpha}[H(x,t)]_{p}-
D_{00}^{\alpha-1}H,$$ 
where 
$$D_{00}^{\alpha-1}H=\lim\limits_{t\to 0}D_{0t}^{\alpha-1}H(x,t),$$
one can get 
\begin{equation}\label{eq23}
[O_{\pm}^{\alpha}H(x,t)]_{k,p}=(p^{\alpha}\pm ik)[H(x,t)]_{k,p}-
[(D_{00}^{\alpha-1}H)(x)]_{k}.
\end{equation}
Applying the double Fourier--Laplace transform to Eq. (\ref{eq35}) and keeping in mind Eq. (\ref{eq23}), 
one can obtain 
$$\left[(O_{+}^{\alpha}-i\gamma	)(O_{-}^{\alpha}-i\gamma	)E_0(x,t)\right]_{k,p}=
(p^{\alpha}-i\gamma	+ ik)[(O_{-}^{\alpha}-i\gamma	)E_0(x,t)]_{k,p}-$$
$$
-\left\{D_{00}^{\alpha-1}[(O_{-}^{\alpha}-i\gamma	)E_0(x,t)]\right\}_k=
(p^{\alpha}-i\gamma	+ ik)(p^{\alpha}-i\gamma	- ik)[E_0(x,t)]_{k,p}-$$
\begin{equation} \label{eq36}
-\Big\{(p^{\alpha}-i\gamma	+ ik)
[D_{00}^{\alpha-1}E_0(x,t)]+
i\sigma D_{00}^{\alpha-1}[ e^{if(x,t)}E_h(x,t)]\Big\}_{k},
\end{equation}
and
$$\left[(O_{-}^{\alpha}-i\gamma	)(O_{+}^{\alpha}-i\gamma	)E_h(x,t)\right]_{k,p}=
(p^{\alpha}-i\gamma	- ik)(p^{\alpha}-i\gamma	+ ik)[E_h(x,t)]_{k,p}-$$
\begin{equation} \label{eq37}
-\Big\{(p^{\alpha}-i\gamma	- ik)D_{00}^{\alpha-1}E_h(x,t)+
i\sigma D_{00}^{\alpha-1}[ e^{-if(x,t)}E_0(x,t)]\Big\}_{k}.
\end{equation}
Here, one has used the relationships ({\it cf.} (\ref{eq31}))
$$D_{00}^{\alpha-1}\left[(O_{-}^{\alpha}-i\gamma	)E_0(x,t)\right] =i\sigma D_{00}^{\alpha-1} [e^{if(x,t)}E_h(x,t)],$$
$$D_{00}^{\alpha-1}\left[(O_{+}^{\alpha}-i\gamma	)E_h(x,t)\right]=i\sigma D_{00}^{\alpha-1}[e^{-if(x,t)}E_0(x,t)].$$
Thus from Eqs. (\ref{eq35}) -- (\ref{eq37}) it directly follows
$$\big({\bf E}(x,t)\big)_{k,p}=\frac{i\sigma}{W}
\Big\{(D_{0t}^{\alpha}
+if'_{x})({\bf KE})-{\bf K}D_{0t}^{\alpha}{\bf E}\Big\}_{k,p}+$$
\begin{equation}\label{eq38}
+\frac{1}{W}\left\{
\left (
\begin{array}{cc}
W^+		& 0  	 \\
0     	& W^-  
\end{array}
\right )
\big(D_{00}^{\alpha-1}{\bf E}\big)_k+i\sigma
\big(D_{00}^{\alpha-1}{\bf KE}\big)_k\right\},
\end{equation}
where
$$W=(p^{\alpha}-i\gamma	)^2+k^2+\sigma^2,
\quad W^{\pm}=p^{\alpha}-i\gamma \pm ik.$$

In \cite{Mamchuev-Chukhovskii-2025-jms} the inverse double Fourier--Laplace transform
\begin{equation}\label{eq381}
\left(\frac{1}{(p^{\alpha}-i\gamma	)^2+k^2+\sigma^2}\right)_{x,t}=
G_{\alpha, \gamma	}(x,t)
\end{equation}
has been taking out for the Green function \cite{Mamchuev-2017-fcaa}
$$G_{\alpha, \gamma}(x,t)=\frac{1}{2}\int\limits_{\lvert x\rvert}^{\infty}
e^{i\gamma	\tau}J_0\left(\sigma\sqrt{\tau^2-x^2}\right) 
\frac{1}{t}\phi\left(-\alpha, 0;-\frac{\tau}{t^{\alpha}}\right)d\tau,$$
here $\phi(\beta, \rho;z)$ is the Wright function
\cite{Wright-1934}
\begin{equation*}
\phi(\beta, \rho;z)=\sum\limits_{k=0}^{\infty}\frac{z^k}{k!\Gamma(\beta k+\rho)},
\quad \beta>-1, \quad \rho\in {\mathbb C}.
\end{equation*}

Let us denote the first and second terms on the right-hand side of relation (\ref{eq38}) as $(({\bf A E})(x,t))_{k,p}$ and $(({\bf B E})(x,t))_{k,p}$, respectively, and then find the inverse double Fourier--Laplace transform of the function 
\begin{equation}\label{eq39}
{\bf E}(x,t)=({\bf A E})(x,t)+({\bf B E})(x,t).
\end{equation}
Using (\ref{eq381}), the integration-by-parts formula, and its analogue in fractional calculus, we obtain for the first term
$$({\bf A E})(x,t)=i\sigma\int\limits_0^t du
\int\limits_{-\infty}^{\infty}G_{\alpha,\gamma}(x-u , t-v)
\Big\{(D_{0v}^{\alpha}
+if'_{u})({\bf KE})-{\bf K}D_{0v}^{\alpha}{\bf E}\Big\}du=$$
$$=i\sigma \int\limits_0^t du \int\limits_{-\infty}^{\infty}
G_{\alpha,\gamma}(x-u , t-v)\frac{\partial}{\partial v}
D_{0v}^{\alpha-1}\Big({\bf K}(u,v){\bf E}(u,v)\Big)du+I_1(x,t)-$$
$$-i\sigma \int\limits_0^t du \int\limits_{-\infty}^{\infty}
G_{\alpha,\gamma}(x-u , t-v){\bf K}(u,v)\frac{\partial}{\partial v}
D_{0v}^{\alpha-1}{\bf E}(u,v)du=$$
$$=i\sigma \int\limits_{-\infty}^{\infty}
\left[G_{\alpha,\gamma}(x-u , t-v)D_{0v}^{\alpha-1}\big({\bf K}{\bf E}\big)(u,v)\right]\Big|_{ v =0}^{ v=t}du-$$
$$-i\sigma \int\limits_0^t du \int\limits_{-\infty}^{\infty}
\frac{\partial}{\partial v}G_{\alpha,\gamma}(x-u , t-v)
\cdot D_{0v}^{\alpha-1}\Big({\bf K}(u,v){\bf E}(u,v)\Big)du+I_1(x,t)-$$
$$-i\sigma \int\limits_{-\infty}^{\infty}
\Big[G_{\alpha,\gamma}(x-u , t-v){\bf K}(u,v)
D_{0v}^{\alpha-1}{\bf E}(u,v)\Big]\Big|_{ v =0}^{ v=t}du+$$
$$+i\sigma \int\limits_0^t du \int\limits_{-\infty}^{\infty}
\frac{\partial}{\partial v}\big[G_{\alpha,\gamma}(x-u , t-v){\bf K}(u,v)\big]
\cdot D_{0v}^{\alpha-1}{\bf E}(u,v)du=$$
$$=I_1(x,t)-I_2(x,t)+I_3(x,t)-$$
$$-i\sigma \int\limits_0^t du \int\limits_{-\infty}^{\infty}
D_{vt}^{\alpha-1}\frac{\partial}{\partial v}G_{\alpha,\gamma}(x-u , t-v)
\cdot {\bf K}(u,v){\bf E}(u,v)du+$$
$$+i\sigma \int\limits_0^t du \int\limits_{-\infty}^{\infty}
D_{vt}^{\alpha-1}\frac{\partial}{\partial v}\big[G_{\alpha,\gamma}(x-u , t-v){\bf K}(u,v)\big]
\cdot {\bf E}(u,v)du,$$
where 
$$I_1(x,t)=i\sigma \int\limits_0^t du \int\limits_{-\infty}^{\infty}
G_{\alpha,\gamma}(x-u , t-v)if'_u(u,v){\bf K}(u,v){\bf E}(u,v)du,$$

$$I_2(x,t)=i\sigma \int\limits_{-\infty}^{\infty}
G_{\alpha,\gamma}(x-u , t)D_{00}^{\alpha-1}({\bf KE})du,$$

$$I_3(x,t)=i\sigma \int\limits_{-\infty}^{\infty}
G_{\alpha,\gamma}(x-u , t)K(u,0)D_{00}^{\alpha-1}{\bf E}du.$$

Ultimately, we obtain
$$({\bf A E})(x,t)=I_1(x,t)-I_2(x,t)+I_3(x,t)+
$$
\begin{equation}\label{eq312}
+i\sigma \int\limits_0^t du \int\limits_{-\infty}^{\infty}
D_{vt}^{\alpha-1}\left[G_{\alpha,\gamma}(x-u , t-v)\frac{\partial}{\partial v}{\bf K}(u,v)\right]
\cdot {\bf E}(u,v)du,
\end{equation}

Taking into account following relations
$$\frac{W^{\pm}}{W}=\frac{q(p)}{q^2(p)+k^2+\sigma^2}\pm \frac{ik}{W}, \quad q(p)=p^{\alpha}-i\gamma,$$ 

$$\left(\frac{p\pm ik}{p^2+k^2+\sigma^2}\right)_{x,t}=
\frac{1}{2}
\left(\frac{\partial}{\partial t}\pm \frac{\partial}{\partial x}\right)
J_0\left(\sigma\sqrt{t^2-x^2}\right)\cdot\Theta[t-|x|],$$

$$\left(e^{-q(p)\tau}\right)_t=e^{i\gamma \tau}\frac{1}{t}
\phi(-\alpha,0;-\tau t^{-\alpha}),$$
and using the Efros's theorem we can get

$$\left(
\left (
\begin{array}{cc}
W^+		& 0  	 \\
0     	& W^-  
\end{array}
\right )
\frac{1}{W}
\right)_{x,t}
=G^1_{\alpha, \gamma}(x,t),
$$
where
$$G^1_{\alpha, \gamma}(x,t)=\frac{1}{2}\int\limits_{\lvert x\rvert}^{\infty}
\frac{1}{t}\phi\left(-\alpha, 0;-\frac{\tau}{t^{\alpha}}\right)
\left(\frac{\partial}{\partial \tau}+\Lambda \frac{\partial}{\partial x}\right)
e^{i\gamma	\tau}J_0\left(\sigma\sqrt{\tau^2-x^2}\right) 
d\tau=$$
$$=i\gamma\cdot G_{\alpha, \gamma}(x,t)-
\frac{1}{2}\int\limits_{\lvert x\rvert}^{\infty}
e^{i\gamma	\tau}
\frac{1}{t}\phi\left(-\alpha, 0;-\frac{\tau}{t^{\alpha}}\right)
\left (
\begin{array}{cc}
\tau+x	& 0  	 \\
0     	&  \tau-x 
\end{array}
\right )
\frac{J_1\left(\sigma\sqrt{\tau^2-x^2}\right)}{\sqrt{\tau^2-x^2}} 
d\tau.$$


From the last relation and (\ref{eq381}) we obtain
\begin{equation}\label{eq313}
({\bf B E})(x,t)=\int\limits_{-\infty}^{\infty}
G^1_{\alpha,\gamma}(x-u , t)D_{00}^{\alpha-1}{\bf E}du+I_2(x,t).
\end{equation}

Taking into account Eqs. (\ref{eq312}), (\ref{eq313}), the  matrix integral equation (\ref{eq39}) may  be  reduce to 
\begin{equation}\label{eq314}
{\bf E}(x,t)-i\sigma\int\limits_{0}^{t}dv\int\limits_{-\infty}^{\infty}
{\bf K}_1(x,t;u,v){\bf E}(u,v)du=
{\bf F}(x,t),
\end{equation}
where 
$$
{\bf K}_1(x,t;u,v)=
D_{vt}^{\alpha-1}\left[G_{\alpha, \gamma}(x-u,t-v)\frac{\partial}{\partial v}{\bf K}(u,v)\right]+$$
$$+if'_{u}(u,v)\cdot G_{\alpha, \gamma}(x-u,t-v){\bf K}(u,v),$$
and the vector ${\bf F}(x,t)$ is equal to
\begin{equation}\label{eq209}
{\bf F}(x,t)=
i\sigma  \int\limits_{-\infty}^{\infty}
\big[G^1_{\alpha,\gamma}(x-u , t)+G_{\alpha, \gamma}(x-u,t){\bf K}(u,0)\big]\varphi(u)du.
\end{equation}

The following theorem holds.

{\bf Theorem 1.} {\it Let $f(x,t),$  $f_x(x,t)$ and $f_t(x,t)$ be the continuous, bounded functions,
$\varphi(x)\in C(-\infty,+\infty),$ and the condition 
$$\lvert\varphi(x)\rvert_*=O\left(\exp\left(\rho\lvert x\rvert^{\varepsilon}\right)\right),
\quad \sigma_0<(1-\alpha)\alpha^{\alpha\varepsilon},
\quad \varepsilon=\frac{1}{1-\alpha};
$$
is satisfied for $\lvert x\rvert\to \infty.$
Then there exists a unique solution to the Cauchy problem (\ref{eq11}), (\ref{eq12}) such that 
\begin{equation}\label{eq321}
D_{0t}^{\alpha}{\bf E}, \frac{\partial {\bf E}}{\partial x}\in C(\Omega), \quad t^{1-\alpha}{\bf E}(x,t)\in C(\overline{\Omega}),
\end{equation}
satisfying the condition
\begin{equation}\label{eq322}
t^{1-\alpha}\lvert{\bf E}(x,t)\rvert_*=O\left(\exp\left(k\lvert x\rvert^{\varepsilon}\right)\right),\quad    as \quad  \lvert x\rvert\to \infty.
\end{equation}
for some $k>0.$
}

In view of the above, the solvability of problem (\ref{eq11}), (\ref{eq12}) is equivalent to the solvability of the system of integral equations (\ref{eq209}).

The estimates \cite{Mamchuev-2017-fcaa}
\begin{equation}\label{eq210}
\left|\frac{\partial^m}{\partial x^m}D_{0t}^{\nu}G_{\alpha, \gamma}(x,t)\right|\leq C \lvert x\rvert^{-\theta}t^{\alpha(1-m+\theta)-\nu-1}, \quad \theta\geq 0, 
\quad t>0,
\end{equation}  
\begin{equation}\label{eq211}
\frac{\partial^m}{\partial x^m}D_{0t}^{\nu}G_{\alpha, \gamma}(x,t)=O\Big( \exp\left(-\sigma_0 \lvert x\rvert^{\varepsilon}t^{-\alpha\varepsilon}\right)\Big), 
\quad \lvert x\rvert\to \infty,
\quad  t<\infty,
\end{equation}  
hold  for arbitrary
$m\in{\mathbb N}\cup \{0\},$ $ \nu\in {\mathbb R}$;
here and further $\sigma_0<(1-\alpha)\alpha^{\alpha\varepsilon},$ $\varepsilon=\frac{1}{1-\alpha};$  $C$ is the positive constant.

We denote by $\lvert A(x,t)\rvert_*$ the scalar function that, for each $(x,t)$, takes the maximum value among all the absolute values of the elements $a_{ij}(x,t)$ of the matrix $A(x,t)$; i.e.,
$\lvert A(x,t)\rvert_* = \max\limits_{i,j}\lvert a_{ij}(x,t)\rvert.$
Likewise, for a vector $B(x, t)$ with components $b_i(x,t)$ we set  
$\lvert B(x,t)\rvert_* = \max\limits_{i}\lvert b_{i}(x,t)\rvert.$

By virtue of the estimates for the function $G_{\alpha, \gamma}(x,t),$
Lemmas 2.1 and 2.2 from the paper \cite{Mamchuev-Chukhovskii-2025-jms}, one can see that
for the kernel $\mathbf{K}_{1}(x,t;u,v)$ of integral equation (\ref{eq314}) and its iterates 
$${\bf K}_n(x,t;\xi,\eta)=\int\limits_{\eta}^{t}dv\int\limits_{-\infty}^{\infty}
{\bf K}_{n-1}(x,t;u,v){\bf K}_1(u,v;\xi,\eta)du.$$
the following estimates 
$$\lvert {\bf K}_n(x,t;\xi,\eta)\rvert_*\leq C
C_{n,\beta,\theta} \lvert x-\xi\rvert^{-\theta}(t-\eta)^{n\beta -1},
$$
$$\lvert {\bf K}_n(x,t;\xi,\eta)\rvert_*\leq C
C_{n,\beta,\theta}(t-\eta)^{n\beta -1}
\exp\left[-\frac{\sigma_0 \lvert x-\xi\rvert^{\varepsilon}}{(t-\eta)^{\alpha\varepsilon}}\right],
$$
hold \cite{Mamchuev-Chukhovskii-2025-jms}, where
$\beta=\alpha\theta,$
$C_{n,\beta,\theta}={\Gamma^n(\beta)C_{1-\theta}^{n-1}}/{\Gamma(n\beta)}.$

These estimates imply the resolvent series convergence
$${\bf R}(x,t;\xi,\eta)=\sum\limits_{n=0}^{\infty}(i\sigma)^n{\bf K}_{n+1}(x,t;\xi,\eta)
$$
estimate for 
\begin{equation}\label{eq317}
\lvert {\bf R}(x,t;\xi,\eta)\rvert_*\leq C (t-\eta)^{\beta-1}
\left(\lvert x-\xi\rvert^{-\theta}+\exp\left[-\frac{\sigma_0 \lvert x-\xi\rvert^{\varepsilon}}{(t-\eta)^{\alpha\varepsilon}}\right]\right).
\end{equation}
and the unique solution existence as well, i.e.,
\begin{equation}\label{eq318}
{\bf E}(x,t)= {\bf F}(x,t)+i\sigma \int\limits_{0}^{t}\int\limits_{-\infty}^{+\infty}{\bf R}(x,t;\xi,\eta) {\bf F}(\xi,\eta)d\xi d\eta
\end{equation}
to equation (\ref{eq314}).
 From the representation of the function $G_{\alpha,\gamma}^{1}$, it is evident that estimates analogous to those in (\ref{eq210}) and (\ref{eq211}) hold for $\lvert G_{\alpha,\gamma}^{1}(x,t)\rvert_*$.

Using the estimates for functions $G_{\alpha,\gamma}$ and $G_{\alpha,\gamma}^{1}$, equation (\ref{eq209}), and lemmas 2.1 and 2.2 from paper \cite{Mamchuev-Chukhovskii-2025-jms}, we obtain
\begin{equation}\label{eq211}
t^{1-\alpha}{\bf F}(x,t)\in C(\overline{\Omega}); \quad
t^{1-\alpha}\lvert{\bf F}(x,t)\rvert_*\leq C \exp\left(k\lvert x\rvert^{\varepsilon}\right). 
\end{equation}
From (\ref{eq317}) and (\ref{eq318}), taking (\ref{eq211}) into account, follow the inclusions (\ref{eq321}) and the relation (\ref{eq322}).

\section{ Cauchy's non-local problem. The case of the FOP $\alpha=1$}
\label{sec:5}
\setcounter{equation}{0}

In the case when FOP  $\alpha=1,$ the basic system (\ref{eq31}) reduces to the following
$$\left (
\begin{array}{cc}
\frac{\partial}{\partial t}-\frac{\partial}{\partial x}	& 0  	 \\
 0   	&      \frac{\partial}{\partial t}+\frac{\partial}{\partial x}
\end{array}
\right )  
\left (
\begin{array}{c}
E_0(x,t)	 \\
E_h(x,t)
\end{array}
\right ) =$$
$$
=i
\left (
\begin{array}{cc}
\gamma						& \sigma \exp[if(x,t)]  	 \\
\sigma\exp[-if(x,t)]    	& \gamma	     
\end{array}
\right )  
\left (
\begin{array}{c}
E_0(x,t)	 \\
E_h(x,t)
\end{array}
\right ).
$$

Note, in limit of the FOP $\alpha=1,$ the function $G_{\alpha, \gamma}(x,t)$ can be cast into the form
$$
G_{1,\gamma}(x,t)=\frac{1}{2}e^{i\gamma t}J_0\left(\sigma \sqrt{t^2-x^2}\right)\Theta[t-\lvert x\rvert]
$$

As was shown earlier \cite{Mamchuev-Chukhovskii-2023-fcaa}, one can find out that for $\alpha\to 1$ integral equation (\ref{eq39}) can be written down into the form
$$
{\bf E}(x,t)=
i\sigma
\int\limits_0^t dv
\int\limits_{x-t+v}^{x+t-v}
{\bf K}_1^1(x, t; u, v)
\cdot {\bf E}(u,v)du+{\bf F}(x,t),$$
where
$${\bf K}_1^1(x, t; u, v)=
\frac{1}{2}e^{i\gamma t}J_0\left(\sigma \sqrt{t^2-x^2}\right)
\left(\frac{\partial}{\partial v}+\Lambda \frac{\partial}{\partial u} 
\right){\bf K}(u,v),$$
$$
{\bf F}(x,t)=\int\limits_{x-t}^{x+t}
\left(\frac{\partial}{\partial t}+\Lambda \frac{\partial}{\partial u} 
+i\sigma {\bf K}(u,0)\right)
\frac{1}{2}e^{i\gamma t}J_0\left(\sigma \sqrt{t^2-x^2}\right)
\cdot{\bf E}(u,0)du.$$

\section{Cauchy's problem  for the displacement field function $f({\bf R})=a x$}
\label{sec:6}
\setcounter{equation}{0}

In this Section, we will build up a solution of the fractional  Cauchy problem with the crystal displacement field function 
$f({\bf R})=a x$.
In this case, we obtain an explicit solution of the problem under study in terms of special functions.

After trivial exponential substitutions  
$$E_0(x,t)=\exp\left(i\frac{a x}{2}\right){\mathcal E}_0(x,t),
\quad
E_h(x,t)=\exp\left(-i\frac{a x}{2}\right){\mathcal E}_h(x,t),$$
the system (\ref{eq11}) reduced to the following system with constant coefficients
\begin{equation}\label{eq61}
\begin{split}
\left (D_{0t}^{\alpha}-\frac{\partial}{\partial x}\right)
{\mathcal E}_0(x,t)=i\left(\gamma+\frac{a}{2}\right)	 {\mathcal E}_0(x,t)+i\sigma {\mathcal E}_h(x,t),  \\
\left (D_{0t}^{\alpha}+\frac{\partial}{\partial x}
\right ){\mathcal E}_h(x,t)=i\sigma {\mathcal E}_0(x,t)+i\left(\gamma+\frac{a}{2}\right) {\mathcal E}_h(x,t).
\end{split}
\end{equation}

According to the works \cite{Mamchuev-2010-du}, \cite{Mamchuev-2012-du} the basic solution of the RL Cauchy problem (\ref{eq11}), (\ref{eq12})
takes the form 
\begin{equation}\label{eq68}
{\bf E}(x,t)=\int\limits_{-\infty}^{\infty}{\bf G}(x,\xi,t)\varphi(\xi)d\xi,
\end{equation}
where the following notations are introduced
\begin{equation}\label{eq69}
{\bf G}(x,\xi,t)=\frac{1}{2}\int\limits_{\lvert x-\xi\rvert}^{\infty}
g(t,\tau){\bf S}(x,\xi,\tau)d\tau+{\bf G}_0(x,\xi,t),
\end{equation}

$$
{\bf S}(x,\xi,\tau)=
\left [\!\!
\begin{array}{cc}
-\sigma e^{iaX_1/2}\frac{\tau-X_1}{\sqrt{\tau^2-X_1^2}}h_1(\tau, X_1)	& 
i\sigma e^{iaX_2/2}h_0(\tau, X_1) \\
i\sigma e^{-iaX_2/2}h_0(\tau, X_1) &			
-\sigma e^{-iaX_1/2}\frac{\tau+X_1}{\sqrt{\tau^2-X_1^2}}h_1(\tau, X_1)
\end{array}
\!\!\right ], 
$$

$$h_0(\tau, X_1)=J_0(\sigma\sqrt{\tau^2-X_1^2}),
\quad 
h_1(\tau, X_1)=J_1(\sigma\sqrt{\tau^2-X_1^2}),$$
$$g(t,\tau)=
\frac{e^{i(\gamma+a/2)\tau}}{t}\phi\left(-\alpha,0; -\tau t^{-\alpha}\right),$$
\begin{equation}\label{eq610}
{\bf G}_0(x,\xi,t)=\left [\!\!
\begin{array}{cc}
e^{-i\gamma X_1}\Theta(-X_1)	& 0 \\
0							& e^{i\gamma X_1}\Theta(X_1)
\end{array}
\!\!\right ]
\frac{1}{t}\phi\left(-\alpha,0; -\lvert X_1\rvert t^{-\alpha}\right), 
\end{equation}
$$X_1=x-\xi, \quad X_2=x+\xi.$$
$\Theta(x)$ is the Heaviside function,
$J_m(z)$ is the $m$-order Bessel function of the argument $z$, $m=0,1$. 

In the article \cite{Mamchuev-Chukhovskii-2025-jms}  
the basic solution of the Gerasimov-Caputo (GC) Cauchy problem been obtained in form
\begin{equation}\label{eq611}
{\bf E}(x,t)=\int\limits_{-\infty}^{\infty}D_{0t}^{\alpha-1}{\bf G}(x,\xi,t)\psi(\xi)d\xi,
\end{equation}
where 
$${\bf E}(x,0)=\psi(x)$$
is the GC Cauchy problem conditions.

\section{Cauchy's non-local problem in the case of the displacement function $f({\bf R})=a x$ and constant initial conditions}

\label{sec:7}
\setcounter{equation}{0}

\subsection{The generalized Euler's formula}

The Stankovic's transformation for the function $v(t)$ is defined as the integral \cite{pskhu-mon-2005}
$$A^{\alpha,\beta}v(t)\equiv
\int\limits_{0}^{\infty}v(\tau)t^{\beta-1}\phi\left(-\alpha,\beta; -\tau t^{-\alpha}\right)d\tau.$$
The following formulas (see \cite[p. 84]{pskhu-mon-2005} for details)
\begin{equation}\label{eq217}
A^{\alpha,\beta}\exp(\lambda t)
= t^{\alpha+\beta-1}E_{\alpha, \alpha+\beta}\left(\lambda t^{\alpha}\right),
\end{equation}
\begin{equation}\label{eq218}
A^{\alpha,\beta}\cos(\lambda t)
= t^{\alpha+\beta-1}E_{2\alpha, \alpha+\beta}\left(-\lambda^2 t^{2\alpha}\right),
\end{equation}

\begin{equation}\label{eq219}
A^{\alpha,\beta}\sin(\lambda t)
= \lambda t^{2\alpha+\beta-1}E_{2\alpha, 2\alpha+\beta}\left(-\lambda^2 t^{2\alpha}\right),
\end{equation}
take place  for any
$\lambda\in {\mathbb C},$
$\alpha\in(0,1),$
$\beta\in {\mathbb R},$
where 
$$E_{\rho, \mu}(z)=\sum\limits_{k=0}^{\infty}\frac{z^k}{\Gamma(\mu+\rho k)}$$
is the Mittag--Leffler-type function
(see \cite[p. 117]{Djrb-1966}).

Using the Eqs. (\ref{eq217}) -- (\ref{eq219}), and the Euler's formula, one obtains the next generalized Euler's formula
$$t^{\alpha+\beta-1}E_{\alpha, \alpha+\beta}\left(i\sigma t^{\alpha}\right)=
t^{\alpha+\beta-1}E_{2\alpha, \alpha+\beta}\left(-\sigma^2 t^{2\alpha}\right)+$$
\begin{equation}\label{eqef}
+i\sigma t^{2\alpha+\beta-1}E_{2\alpha, 2\alpha+\beta}\left(-\sigma^2 t^{2\alpha}\right),
\end{equation}
which will be useful below.
Denoting the parameter $\mu = \alpha + \beta,$ equality (\ref{eqef}) can be rewritten as follows
\begin{equation}\label{eqef1}
t^{\mu-1}E_{\alpha, \mu}\left(i\sigma t^{\alpha}\right)=
t^{\mu-1}E_{2\alpha, \mu}\left(-\sigma^2 t^{2\alpha}\right)+
i\sigma t^{\alpha+\mu-1}E_{2\alpha, \alpha+\mu}\left(-\sigma^2 t^{2\alpha}\right).
\end{equation}
Note that for $\mu=1$, formula (\ref{eqef1}) was obtained in \cite{Gutov-2005}.

Introducing the following notations
$$
exp_{\alpha,\beta}(\sigma, t)=t^{\alpha+\beta-1}
E_{\alpha, \alpha+\beta}\left(\sigma t^{\alpha}\right), $$
$$\cos_{\alpha,\beta}(\sigma, t)=t^{\alpha+\beta-1}
E_{2\alpha, \alpha+\beta}\left(-\sigma^2 t^{2\alpha}\right), $$

$$\sin_{\alpha,\beta}(\sigma, t)=\sigma t^{2\alpha+\beta-1}
E_{2\alpha, 2\alpha+\beta}\left(-\sigma^2 t^{2\alpha}\right),
$$
one rewrites formula (\ref{eqef}) into the form 
$$exp_{\alpha,\beta}(i\sigma, t)=\cos_{\alpha,\beta}(\sigma, t)+
i \sin_{\alpha,\beta}(\sigma, t),$$
$$exp_{\alpha,\beta}(-i\sigma, t)=\cos_{\alpha,\beta}(\sigma, t)-
i \sin_{\alpha,\beta}(\sigma, t).$$

From the last two equalities follow the next relations
\begin{equation}\label{eq71gcos}
\cos_{\alpha,\beta}(\sigma, t)=
\frac{exp_{\alpha,\beta}(i\sigma, t)+exp_{\alpha,\beta}(-i\sigma, t)}{2},
\end{equation}
\begin{equation}\label{eq71gsin}
\sin_{\alpha,\beta}(\sigma, t)=
\frac{exp_{\alpha,\beta}(i\sigma, t)-exp_{\alpha,\beta}(-i\sigma, t)}{2i}.
\end{equation}

\subsection{Constant initial conditions }

Let us consider the case when 
{the initial conditions for the wavefield amplitudes} 
are constant, and they are equal to
\begin{equation}\label{eq71}
\lim\limits_{t\to 0}D_{0t}^{\alpha-1}E_0(x,t)=\varphi_0(x)\equiv C_0, \quad
\lim\limits_{t\to 0}D_{0t}^{\alpha-1}E_h(x,t)=\varphi_h(x)\equiv C_h.
\end{equation}
Accordingly, formulae (\ref{eq68})--(\ref{eq610}) can be simplified and expressed in the terms of the Mittag--Leffler-type functions.

Bearing in mind Eqs. (\ref{eq68}), (\ref{eq71})
and the work \cite{Mamchuev-Chukhovskii-2025-jms} calculations details, the vector ${\bf E}(x, t)$ can be cast in the form
$$
{\bf E}(x,t)=\int\limits_{-\infty}^{\infty} 
{\bf G}(x,\xi,t)\left (
\begin{array}{c}
C_0	 \\
C_h	
\end{array}
\right )d\xi=$$
\begin{equation}\label{eq7100}
=\int\limits_0^{\infty} 
\frac{1}{t} \phi\left(-\alpha,0;-\frac{\eta}{t^{\alpha}}\right)
\Big[
{\bf N}_1(x)e^{i(a_1+k)\eta}
+{\bf N}_2(x)e^{i(a_1-k)\eta}\Big]
\left (
\begin{array}{c}
C_0	 \\
C_h	
\end{array}
\right )d\eta
\end{equation}
where
$${\bf N}_{1,2}(x)=\frac{1}{2}\left (
\begin{array}{ll}
1\mp \frac{a}{2k}&  
\pm\frac{\sigma}{k} e^{iax} \\
\pm\frac{\sigma}{k} e^{-iax}	&  
1\mp\frac{a}{2k}
\end{array}
\right ), \quad
a_1=\gamma+\frac{a}{2},  \quad k=\sqrt{a^2/4+\sigma^2}.$$
Thus, from Eq. (\ref{eq7100}) for the wave vector ${\bf E}(x,t)$ one has obtained
\begin{equation}\label{eq7101}
{\bf E}(x,t)=t^{\alpha-1}\Big[
{\bf N}_1(x)E_{\alpha, \alpha}\left(i(a_1+k)t^{\alpha}\right)
+{\bf N}_2(x)E_{\alpha, \alpha}\left(i(a_1-k)t^{\alpha}\right)\Big]
\left (
\begin{array}{c}
C_0	 \\
C_h	
\end{array}
\right ),
\end{equation}
and or in the aggregate form 
\begin{equation}\label{eq7103}
{\bf E}(x,t)=
\overline{{\bf N}}(x,t)
\left (
\begin{array}{c}
C_0	 \\
C_h	
\end{array}
\right )
\end{equation}
where
\begin{equation}\label{eq71031}
\overline{{\bf N}}(x,t)=\frac{1}{2}\left (
\begin{array}{lr}
(E^++E^-)-\frac{a}{2k}(E^+-E^-)&  
\frac{\sigma}{k} e^{iax}(E^+-E^-) \\
\frac{\sigma}{k} e^{-iax}(E^+-E^-)	&  
(E^++E^-)-\frac{a}{2k}(E^+-E^-)
\end{array}
\right ),
\end{equation}
$$E^+=exp_{\alpha,0}(i(a_1+k),t), \quad
E^-=exp_{\alpha,0}(i(a_1-k),t).$$

\subsection{The two particular cases of the RL Cauchy problem solution}

In the special case when
$$
\lim\limits_{t\to 0}D_{0t}^{\alpha-1}E_0(x,t)=1, \quad
\lim\limits_{t\to 0}D_{0t}^{\alpha-1}E_h(x,t)=1,
$$
formula (\ref{eq7101}) takes the form

\begin{equation}\label{eq7131}
{\bf E}(x,t)=\frac{t^{\alpha-1}}{2}
\left (
\begin{array}{ll}
1-\frac{a}{2k}+\frac{\sigma}{k} e^{iax}	&  
1+\frac{a}{2k}-\frac{\sigma}{k} e^{iax}\\
1-\frac{a}{2k}+\frac{\sigma}{k} e^{-iax}	&  
1+\frac{a}{2k}-\frac{\sigma}{k} e^{-iax}
\end{array}
\right )
\left (
\begin{array}{c}
E_{\alpha, \alpha}\left(i(a_1+k)t^{\alpha}\right)	 \\
E_{\alpha, \alpha}\left(i(a_1-k)t^{\alpha}\right)
\end{array}
\right ).
\end{equation}

Accordingly, for the following conditions as 
$$
\lim\limits_{t\to 0}D_{0t}^{\alpha-1}E_0(x,t)=1, \quad
\lim\limits_{t\to 0}D_{0t}^{\alpha-1}E_h(x,t)=0,
$$
from Eq. (\ref{eq7101}), one has obtained 
\begin{equation}\label{eq713}
{\bf E}(x,t)=\frac{t^{\alpha-1}}{2}
\left (
\begin{array}{lr}
1-\frac{a}{2k}	&  
1+\frac{a}{2k}\\
\frac{\sigma}{k} e^{-iax}	&  
-\frac{\sigma}{k} e^{-iax}
\end{array}
\right )
\left (
\begin{array}{c}
E_{\alpha, \alpha}\left(i(a_1+k)t^{\alpha}\right)	 \\
E_{\alpha, \alpha}\left(i(a_1-k)t^{\alpha}\right)
\end{array}
\right ).
\end{equation}

It is interesting that in the case of FOP $\alpha=1$ 
from (\ref{eq71031}) one can see that 
the matrix $\overline{{\bf N}}(x,t)$ takes the form
$$\overline{{\bf N}}(x,t)={\bf N}(x,t)=e^{ia_1t}\left (
\begin{array}{lr}
\cos kt -\frac{a}{2k}i\sin kt&  
\frac{\sigma}{k} e^{iax}i\sin kt \\
\frac{\sigma}{k} e^{-iax}i\sin kt	&  
\cos kt-\frac{a}{2k}i\sin kt
\end{array}
\right ).$$

If in addition $a=0$ the vector ${\bf E}(x,t)$ is equal to
\begin{equation}\label{eq7103}
{\bf E}(x,t)=e^{i\gamma t}
\left (
\begin{array}{lr}
\cos kt &  
i\sin kt \\
i\sin kt	&  
\cos kt
\end{array}
\right )
\left (
\begin{array}{c}
C_0	 \\
C_h	
\end{array}
\right )=
e^{i\gamma t}
\left (
\begin{array}{lr}
C_0&  
C_h \\
C_h	&  
C_0
\end{array}
\right )
\left (
\begin{array}{c}
\cos kt	 \\
i\sin kt	
\end{array}
\right ).
\end{equation}

In the case when the FOP $\alpha\to 1$, the solution (\ref{eq713}) has reduced to
\begin{equation}\label{eq714}
{\bf E}(x,t)=e^{ia_1t}
\left (
\begin{array}{l}
\cos{kt}-\frac{a}{2k}i\sin{kt}\\
\frac{\sigma}{k} e^{-iax}i\sin{kt}
\end{array}
\right ).
\end{equation}
Finally, putting on the parameters $\gamma=0$ and $a=0,$  the solution (\ref{eq714}) goes to
'pendulum solutions' of the canonic TT equations regarding the Bloch wave-field
amplitudes in the perfect crystal (cf.  \cite{Mamchuev-Chukhovskii-2023-fcaa})
$$E_0=\cos \sigma t, \quad E_h=i\sin \sigma t.$$

In the case when $a=\gamma=0$ from (\ref{eq71031}) taking into account formulas (\ref{eq71gcos}) and (\ref{eq71gsin}) one obtains the following
"quasi pendulum solutions": for the RL Cauchy nonlocal problem
\begin{equation}\label{eq710333}
{\bf E}(x,t)=
\left (
\begin{array}{lr}
C_0&  
C_h \\
C_h	&  
C_0
\end{array}
\right )
\left (
\begin{array}{c}
\cos_{\alpha,0} (k,t)	 \\
i\sin_{\alpha,0} (k,t)	
\end{array}
\right ),
\end{equation}
and for GC Cauchy problem
\begin{equation}\label{eq710333}
{\bf E}(x,t)=
\left (
\begin{array}{lr}
C_0&  
C_h \\
C_h	&  
C_0
\end{array}
\right )
\left (
\begin{array}{c}
\cos_{\alpha,1-\alpha} (k,t)	 \\
i\sin_{\alpha,1-\alpha} (k,t)	
\end{array}
\right )
\end{equation}
respectively.

\section{Physical interpretation}

Here, an important remark has to be done regarding the fundamental 
differences between the GC and RL approaches in consideration to solving the Cauchy problem in the fractional X-ray diffraction optics.

The use of the Gerasimov -- Caputo fractional derivative makes it possible to take into account the effect of spatial nonlocality (modeled by the power-law function $t^{-\alpha}/\Gamma(1-\alpha)$) on the alternation  rate of the amplitudes of the transmitted and diffracted waves along the crystal depth $t$, i.e., the effect of nonlocality on the gradient of evolution of these amplitudes. Therefore, the initial values 
${\bf E}(x,0)$ of these amplitudes themselves are specified.

The Riemann -- Liouville fractional derivative describes the evolution of the wave field along the coordinate $t$  as the alternation rate of  the amplitudes under the assumption that they are affected by spatial nonlocality (also  modeled by a power-law function), i.e., the alternation rate of the "nonlocal" amplitudes. Therefore, the initial values of the "nonlocal" amplitudes are specified, i.e., the initial value of the fractional integral $\lim\limits_{t\to 0}D_{0t}^{\alpha-1}{\bf E}(x, t)$. In order for this initial value to be specified correctly, the function ${\bf E}(x, t)$ must have a power-law singularity at zero whose order is exactly $1-\alpha$.

Thus, in both the GC and RL approaches under consideration, the corresponding terms in the fractional-derivative-based Takagi-Taupin (TT) equations carry the different physical sense, namely: 
1) in the GC-approach case, they represent the "nonlocal" alternation rate of the wavefield amplitudes,
2) in the RL-approach, they represent the alternation rate of the "nonlocal" wavefield amplitudes. 
The above means that they each from other have the different physics sense. Accordingly, the first ones can be termed as the "local" (in the GC-approach) and the second ones as "nonlocal" (in the RL-approach) wavefield amplitudes. 

Hereby, it should also be noted that in the limiting case $\alpha=1$, both the GC and RL approaches yield the so-called standard solution of the Cauchy problem in the X-ray diffraction optics.

\section{Conclusion}

In the framework of the RL derivative concept, based on the resolvent-solution of the integral Fredholm -- Volterra matrix equation, a new theoretical approach of the processing the fractional X-ray diffraction crystal optics data has been elaborated.
Noteworthy is the fact that the worldwide X-ray diffraction crystal optics related to the fractional order parameter (FOP) $\alpha$ equal to unity, has been  reformulated using the fractional derivatives concept  
in the Riemann -- Liouville sense $D_{0t}^{\alpha}{\bf E}$ with an arbitrary FOP $\alpha$, $\alpha \in(0, 1]$. 
 
On the example of the X-ray diffraction by perfect crystal, the function $f(r) =ax$, one has found the exact analytical matrix solution of the RL-approach-Cauchy problem for an arbitrary value of the FOP 
$\alpha$, $\alpha \in(0, 1]$.
In the case of constant initial data, these solutions are expressed in terms of Mittag-Leffler-type functions.
The latter have been analyzed in comparing with the corresponding productive solutions obtained in the Gerasimov -- Caputo fractional calculus, respectively. 

In concluding, as a viewpoint of fundamental physics as well an endeavor of revealing the non-locality of the X-ray quanta-medium atoms interaction feature, authors could suggest a new fractional integral-differential derivative of kind
$$D_{at}^{\alpha,\lambda} g(t)=\frac{1}{\Gamma(-\alpha+1)}\frac{d}{dt} \int\limits_a^t\frac{g(s)ds}{(t-s+\lambda)^{\alpha}},  \quad 
0<\alpha\leq 1, \quad
0\leq\lambda\leq 1.$$

In the case when $\lambda=0$, the derivative $D_{at}^{\alpha,\lambda} g(t)$ goes to the RL one $D_{at}^{\alpha} g(t)$ exactly.
Clearly, using of the above fractional integral-differential derivative $D_{at}^{\alpha,\lambda} g(t)$ is of a definite interest for solving the non-local X-ray diffraction Cauchy problem being nothing else a good topic for future job.

\end{document}